\documentclass{amsart}
\title [Weak mixing of maps with bounded cutting parameter]{Weak mixing of maps with bounded cutting parameter} % If the title is too long for the running heads, give a
\author[E.H. El Abdalaoui]{El Houcein EL ABDALAOUI}   % A shortened version for the running heads may be given in []'s,
\address{Laboratoire de Math\'ematiques Rapha\"el Salem \\
UMR 6085 CNRS -- Université de Rouen\\
site Colbert \\
76821 Mont-Saint-Aignan  CEDEX\\
France}
\email{ElHocein.Elabdalaoui@univ-rouen.fr}   % You may include the URL for your home page along with
\author[A. Nogueira]{Arnaldo NOGUEIRA${}^{*}$}  % A shortened version for the running heads may be given in []'s,
\address{Instituto de Mathem\'atica, Universidade Federal do Rio de Janeiro\\
Caixa Postal 68530\\
 21945-970 Rio de Janeiro, RJ- Brazil.}
\email{nogueira@acd.ufrj.br}
\address{IML - Marseille\\
163, av. de Luminy - Case 907,\\
 13288 Marseille Cedex 9, FRANCE
}
\email{nogueira@iml.univ-mrs.fr}   % You may include the URL for your home page along with
\author[T. de la Rue]{Thierry DE LA RUE}   % A shortened version for the running heads may be given in []'s,
\address{Laboratoire de Math\'ematiques Rapha\"el Salem \\
UMR 6085 CNRS -- Université de Rouen\\
site Colbert \\
76821 Mont-Saint-Aignan  CEDEX\\
France}
\email{thierry.de-la-rue@univ-rouen.fr}   % You may include the URL for your home page along with
\thanks {(*) Partially supported by grants from CNPq-Brazil 301456/80
and FINEP/ CNPq/MCT 41.96.0923.00 (PRONEX).}
\keywords{Weak mixing, Ornstein transformations, interval exchange
maps}

\subjclass[2000]{ Primary : 37D05; secondary : 37A30.}
\newtheorem{thm}{Theorem}[section]

\newtheorem{lem}[thm]{Lemma}
\newtheorem{coro}[thm]{Corollary}

\theoremstyle{definition}

\def\egdef{\stackrel{\rm def}{=}}
\def\tower{
\setlength{\unitlength}{1mm}
\begin{picture}(0,0)
\put (0,0) {\framebox (50,20)}
\put (3,-4) {\makebox (3,-4)[br] {$B_{k}$}}
\put (16,-6) {\makebox (16,-6)[b]{$\underbrace {
\quad  \quad   \quad  \quad  \quad  \quad  \quad  \quad \quad  \quad
\quad  \quad  \quad}_{B_{k-1}}$}}
\put (17,-16) {\makebox (17,-16)[b]{$\overline {
\quad  \quad  \quad \quad ~ ~ p_{k-1}
\quad  \quad   \quad  \quad \quad   \quad ~ ~
}$}}
\multiput (5,0)(5,0){10}{\line (0,20){20}}
\put (0,4) {\line (4,0){50}}
\put (3,0){\vector (0,1){4}}
\thicklines
\put (0,0) {\line (5,0){5}}
\multiput (0,21)(0,1){5}{\line (5,0){5}}
\put (-3,27){\makebox (-2,27)[bl] {$a_1^{(k-1)}$}}
\put (5,21){\line(5,0){5}}
\put (5,22){\line(5,0){5}}
\put (8,24){\makebox (8,24)[bl] {$a_2^{(k-1)}$}}
\put (10,21){\makebox (10,21)[b] {$\cdots \cdots$}}
\multiput (25,21)(0,1){6}{\line (25,0){5}}
\put (23,28){\makebox (23,28)[bl] {$a_j^{(k-1)}$}}
\put (26,21){\makebox (26,21)[b] {$\cdots \cdots$}}
\multiput (45,21)(0,1){3}{\line (45,0){5}}
\put (45,25){\makebox (45,25)[bl] {$a_{p_{k-1}}^{(k-1)}$}}
\put (-40,-25){\centerline{Figure 1: $k^{\hbox{th}}$--tower.}}
\end{picture}}

\def\N{\mathbb{N}}
\def\Z{\mathbb{Z}}
\def \T{\mathbb{T}}
\def \P{\mathbb{T}}

\def\build#1_#2^#3{\mathrel{\mathop{\kern 0pt#1}\limits_{#2}^{#3}}}
\def\tend#1#2{\build\hbox to 12mm{\rightarrowfill}_{#1\rightarrow #2}^{}}
\def\netendpas#1#2{\build\hbox to 12mm{$\not \longrightarrow$}_{#1 \rightarrow
#2}^{}}

\def \converge#1#2#3{\build\hbox to 15mm {\rightarrowfill}_{#1\rightarrow #2}^
{\hbox{\scriptsize #3}}}

\begin{document}
       % Abstracts are required, and will be distributed by a listserv list.
       % Please place the text of your abstract in the next environment.

\begin{abstract}
In the class of Ornstein transformations the mixing property
satisfies a 0-1 law. Here we consider Ornstein's construction
with bounded cutting parameter. In fact, these latter
transformations are not mixing, however it is proved that the weak
mixing property occurs with probability one. Our situation is
similar to the case of interval exchange transformations whose
link with the cutting and stacking
construction relies in a dynamical process called Rauzy induction.
\end{abstract}

\maketitle
%\setcounter{tocdepth}{1}
%\tableofcontents

\section{Introduction}

\noindent{}In a nowadays classical work \cite{Ornstein}, Ornstein
associates to every point $\omega$ in a probability space $\Omega$
a rank one transformation $T_{\omega}$ which he proves to be
mixing for almost every $\omega$. There are many extensions and
generalizations of Ornstein's construction. Here we regard a
class of natural examples of rank one maps of the interval for
inspiration to consider a generalization of Ornstein's result.

It has been established by Veech \cite{Veech} that interval
exchanges are almost surely rank one, assuming the permutation is
irreducible. The link between the cutting and stacking
construction and interval exchanges is done through Rauzy
induction which is a way to define first return induced
transformation without increasing the number of exchanged
intervals (see \cite{Rauzy}). However Katok \cite{Katok2} proved
that no interval exchange map is mixing. Katok's result remains at
the present time, the only universal result about the spectrum of
interval exchanges. In fact, he showed that every interval
exchange map is $\alpha$-rigid (a transformation $T$ is said to be
$\alpha$-rigid, $0 < \alpha <1$, if there exists an infinite
sequence of integers $\{n_k\}_{k \in \N}$ such that $
\lim_{k\rightarrow\infty}\mu(T^{n_k}A \cap A) \geq \alpha \mu(A)$,
for every measurable set $A$). Naturally, one asks the following
question:\\

{\it {Question 1.1.  Does any interval exchange map have singular
spectral type?}}\\

It is known that the stacking construction using a constant
cutting parameter results in a map which is not mixing. In the
class of interval exchanges, whether weak mixing is almost surely
satisfied was still an open question (see \cite{Rauzy}) until the
recent work of Avilla and Forni \cite{Avilla}. We recall that
Katok and Stepin \cite{Katok1} and Veech \cite{Veech} have proved
that for some permutations almost every interval exchange map is
weak mixing, nevertheless in \cite{Nogueira} it is shown that
those permutations force the eigenvalue to be 1.

Here we will consider the class of Ornstein transformations with
bounded cutting
parameter.
We prove that in this case the weak mixing property occurs with
probability one.

It is clear that our result is in the same spirit of the result of
\cite{Avilla}. One may hope that there is sufficient analogy
between our construction and interval exchange maps for our
methods to extend to this case and to obtain a unified proof.

Here we begin by considering a generalization of Ornstein's
construction with less restrictions (namely in Ornstein's case,
some non explicit cutting parameters satisfying some growth
condition ensuring mixing are shown to exist, whereas here we fix
them in advance). In this context, answering a question asked by
J.-P. Thouvenot, we extend a result proved in [8] to the class of
Ornstein transformations with bounded cutting parameter.

We will assume that the reader is familiar with the method of cutting
and stacking for constructing rank one transformations.

\section{Construction of rank one transformations}

\noindent Using the cutting and stacking method described in
Friedman \cite{Friedman1,Friedman2}, we can define inductively a
family of measure preserving rank one transformations, as follows:
\vskip 0.1cm \noindent Let $B_0$ be the unit interval equipped
with the Lebesgue measure. At the first stage $B_0$ is splitted
into $p_0$ equal parts, add spacers and form a stack of height
$h_{1}$ in the usual manner. At the $k$-th stage we divide the
stack obtained at the $(k-1)$-th stage into $p_{k-1}$ equal
columns, add spacers and obtain a new stack of height $h_{k}$. If
during the $k$-th stage of our construction the number of spacers
put above the $j$-th column of the $(k-1)$-th stack is
$a^{(k-1)}_{j} \geq 0$, $1\leq j \leq p_{k}$, then we have
$$h_{k} = p_{k-1}h_{k-1} + \sum_{j=1}^{p_{k-1}}a_{j}^{(k-1)}.$$
\vskip 3 cm
\hskip 3.5cm
\tower
\vskip 3.0cm
%\centerline{Figure 1:$(n+1)^{\hbox{th}}$--tower.}
\noindent{}Proceeding in this way we get a rank one transformation
$T$ on a certain measure space $(X,{\mathcal B} ,\nu)$ which may
be finite or $\sigma-$finite depending on the number of spacers
added.

The construction of a rank one transformation thus needs two
parameters $(p_k)_{k=0}^\infty$ (cutting and stacking parameter) and
$((a_j^{(k)})_{j=1}^{p_k})_{k=0}^\infty$ (spacers parameter). We define
$$T \egdef T_{(p_k, (a_j^{(k)})_{j=1}^{p_k})_{k=0}^\infty}.$$

\section{Ornstein's class of transformations}

In Ornstein's construction, the $p_k$'s are rapidly increasing and
the number of spacers, $a_i^{(k)}$, $1 \leq i\leq p_k-1$, is
chosen stochastically in a certain way (subject to certain
bounds). This may be organized in various ways as noted by
Bourgain \cite{Bourgain}, in fact, let $(t_k)$ be a sequence of
positive integers such that $2t_k \leq h_{k}$. We choose now
independently, using the uniform distribution on the set
$X_k=\{-t_k,\cdots,t_k\} $, the numbers $(x_{k,i})_{i=1}^{p_k-1}$,
and $x_{k,p_k}$ is chosen deterministically in $\N$. We set, for
$1 \leq i \leq p_k$,
$$
a_i^{(k)} =2t_{k} + x_{k,i} - x_{k,i-1}, ~~{\rm with} ~~x_{k,0} =
0.
$$
Then one sees that
$$
h_{k+1} = p_k(h_k +2 t_k) + x_{k,p_k}.
$$
So the deterministic sequences of positive integers
$(p_k)_{k=0}^\infty$ and $(x_{k,p_k})_{k=0}^\infty$ completely determine the
sequence of heights $(h_k)_{k=1}^\infty$. The total measure of the resulting
measure space is finite if $\displaystyle
\sum_{k=1}^{\infty}\frac{t_k}{h_k}+\sum_{k=1}^\infty \frac{x_{k,p_k}}{p_kh_k} <
\infty$. We will assume that this requirement is satisfied.

We thus have a probability space of Ornstein transformations
$\prod_{k=1}^\infty X_k^{p_k-1}$ equip\-ped with the natural
probability measure $\P \egdef\otimes_{k=1}^{\infty} P_k$, where
$P_k\egdef\otimes_{l=1}^{p_k-1}{\mathcal {U}}_k$; ${\mathcal {U}}_l$ is
the uniform probability on $X_l$. We denote this space by
$(\Omega, {\mathcal {A}}, {\P})$. The projection of $\Omega$ onto
the $i$-th co-ordinate space of $\Omega_k \stackrel {\rm def} {=}
X_k^{p_k-1}$, $1 \leq i \leq p_k-1$ is $x_{k,i}$. Naturally each
point $\omega =(\omega_k
=(x_{k,i}(\omega))_{i=1}^{p_k-1})_{k=1}^\infty$ in $\Omega$
defines the spacers and therefore a rank one transformation which
we denote by $T_{\omega,x}$, where $x=(x_{k,p_k})_{k=1}^{\infty}$
is {\it admissible}, i.e.,
\[
\sum_{k=1}^\infty\frac{x_{k,p_k}}{p_kh_k} < \infty.
\]

The above construction gives a more general definition of the random
construction due to Ornstein.

We recall that an automorphism is said to be {\it totally ergodic
} if all its nonzero powers are ergodic. It is shown in
\cite{elabdal} that the classical Ornstein \linebreak[1]
transformations are almost surely totally ergodic using the fact
that a measure preserving automorphism is totally ergodic if and
only if no root of unity other than 1 is an eigenvalue. In fact it
is proved in \cite{elabdal} that, for a fixed $z \in \T
\setminus\{1\} \equiv [0,1)\setminus \{0\}$, $\{\omega:z {\rm { ~
is ~an ~eigenvalue ~of ~}} T_{\omega} \}$ is a null measure set.
Later, in \cite{elabdal1}, by Van der Corput's inequality and
Bernstein's inequality on the derivative of a trigonometric
polynomial combined with the ingredients of \cite{elabdal} and
\cite{Blum}, we obtain a null measure set $N$ such that for all
$\omega \notin N,$ $T_{\omega}$ has no eigenvalue other than 1
provided that for infinitely many $n'$s we have $t_{n} = p_{n}$
($p_n$ goes to $\infty$, as $n$ goes to $\infty$), and this
implies the almost sure weak mixing property.

We note that it is an easy exercise to show that the spectral
properties satisfy the {\it Zero-One law}. We shall denote by
${\mathcal {WMIX}}$, the $\omega$ set for which $T_{\omega}$ is
weak mixing.

\section{Ornstein transformations with bounded cutting parameter}

Here we assume that the cutting parameter $(p_k)_{k \geq 0}$ is
bounded. Next we state our main result.\\

%\noindent{\bf Theorem 4.1.}  {\it
\begin{thm}
Let $x={(x_{k,p})}_{k \in \N}$
be admissible sequences of positive integers (i.e.$\displaystyle
\sum_{k=1}^\infty\frac{x_{k,p_k}}{p_kh_k} < \infty$) Then
$\P({\mathcal {WMIX}})=1$.
\end{thm}

First we remark that it is an easy exercise to see that the weak
mixing property
occurs with probability 1 if the series $\displaystyle
\sum_{k=1}^{\infty}\frac1{t_k}$ diverges.  In fact, one can show that
Chacon's pattern occurs for infinitely many values of $k$ with probability 1.
Hence we assume that the series $\displaystyle \sum_{k=1}^{\infty}\frac1{t_k}$
converges.

A standard argument yields that for any rank one transformation
$T_{(p_k, (a_j^{(k)})_{j=1}^{p_k})_{k=0}^\infty}$ with bounded
cutting parameter, if $\lambda=e^{2i\pi\alpha}$ is a eigenvalue
then
\[
\lambda^{h_k+a_1^{(k)}} \tend{k}{\infty}1 {\rm
{~~or~~}}{||(h_k+a_1^{(k)})\alpha||}\tend{k}{\infty}0,
\]
where $||x||=d(x,\Z)$.

Let $\mathcal{N}$ be the
subsequence of positives integers $(n_k)$ and put
\[
G({\mathcal{N}})\egdef\{\lambda=e^{2i\pi\alpha}  \in \T
~:~||n_k\alpha||\tend{k}{\infty}0\}.
\]

\subsection*{4.1. Some general facts about $G(\mathcal{N})$ when $\displaystyle
\left (\frac{n_{k+1}}{n_k}\right)_{k \in \N} $ is bounded}

Let $\varepsilon$ be a positive number and $L$ a positive integer. Put
\[
A_{(n_k),L}^{(\varepsilon)}=\{\lambda \in [0,1)~:~ ||n_k\lambda||
<\varepsilon, \forall n_k >L\}.
\]
Observe that we have
\begin{eqnarray}{\label {eqn : G1}}
G({\mathcal{N}})=\bigcap_{\epsilon >0}\bigcup_{L \in
\N}A_{(n_k),L}^{(\varepsilon)}.
\end{eqnarray}

\noindent{}Now, from this observation, we shall study the
properties of $G(\mathcal{N})$ when the sequence
$\left(\displaystyle \frac{n_{k+1}}{n_k} \right)$ is bounded. We
note, first, that we have the
following lemma\\

\begin{lem}
Assume that there exists a
postive number $M$ such that
\[
\frac{n_{k+1}}{n_k} < M, {\rm {~for~ any ~}} k \in \N.
\]
 Then, for any $\varepsilon$ less than $\displaystyle \frac{1}{4M}$ we
 have $|A_{(n_k),L}^{(\varepsilon)}| \leq n_{k_0}$, where $k_0$ is the smallest positive integer such that
 $n_{k_0}>L$ and $|A_{(n_k), L}^{(\varepsilon)}|$ is the cardinal of
 $A_{(n_k),L}^{(\varepsilon)}$.
\end{lem}

\begin{proof} Observe that

\begin{eqnarray}{\label {eqn : G2}}
A_{(n_k),L}^{\varepsilon}=\bigcap_{n_k \geq
L}B_{n_k}^{\varepsilon},
\end{eqnarray}

\noindent{}where $B_{n_k}^{\varepsilon}=\{\lambda \in
[0,1)~:~||n_k \lambda||<\varepsilon\}.$ We deduce form the
definition of  $k_0$ and (\ref{eqn : G2}) that
\[
A_{(n_k),L}^{\varepsilon}=\bigcap_{k \geq
k_0}B_{n_k}^{\varepsilon}.
\]
\noindent{}But $B_{n_k}^{\varepsilon}$ is the union of the
intervals centered on $\displaystyle \frac{j}{n_k}$, $0 \leq j
\leq n_k-1$ and of length $\displaystyle
\frac{2\varepsilon}{n_k}$. It follows that if $I$ is some interval
from $B_{n_k}^{\varepsilon}$ centered on some $\displaystyle
\frac{j}{n_k}$ has a non-empty intersection with two different
intervals from $B_{n_{k+1}}^{\varepsilon}$, then we must have

\begin{eqnarray}{\label {eqn : G3}}
\frac{1-2\varepsilon}{n_{k+1}} \leq \frac{2\varepsilon}{n_k},
\end{eqnarray}
It follows from (\ref{eqn : G3}) that
\[
\frac{2\varepsilon}{1-2\varepsilon}\geq \frac1{M} {\rm
{~~hence~~}} 4\varepsilon \geq \frac1{M},
\]
which yields a contradiction. Now, Let $x,x'$ be in
$A_{(n_k),L}^{\varepsilon}$. Assume that
$x,x'$ are in the same interval from $B_{n_{k_0}}^{\varepsilon}$;
Then from the above we deduce by induction that $x,x'$ are in the same interval
from $B_{n_k}^{\varepsilon}$, for any $k \geq k_0$. It follows
that $x=x'$ and this yields that $|A_{(n_k),L}^{(\varepsilon)}|
\leq n_{k_0}$. The proof of the lemma is complete.
\end{proof}

\begin{coro}
If $\displaystyle \left(\frac{n_{k+1}}{n_k} \right)$ is bounded then $G(n_k)$ is
countable.
\end{coro}

\noindent{}We have also the following lemma.

\begin{lem}
\label{4.1.2}
Let ${\mathcal {N}}=\{n_k, k \in \N\}$ and ${\mathcal
{N'}}=\{n'_k, k \in \N\}$ two sequences of positive integers such
that there exist infinitely many $k$ for which : $n'_k=n_k+1$,
then $G({\mathcal {N}}) \bigcap G({\mathcal {N}}') = \{1\}$.
 \end{lem}

\begin{proof}
Straightforward. (Chacon's argument!)
\end{proof}

\section{Application and proof of theorem.}
First, set
\[
n_k(\omega)\egdef h_k+x_{k,1}(\omega), {~~and~~} {\mathcal
{N}}(\omega)\egdef \{n_k(\omega), k\in \N\}~~~~\omega \in \Omega.
\]and observe that we have, for any $\omega \in \Omega$,
\[
\frac{n_{k+1}}{n_k} \leq p+1.
\]Let ${\mathcal{F}}$ be a the $\sigma$-algebra
generated by the random variables $\{{x_{2k}}, k\in \N \}$ and the
event $x_{2k+1}\in A$, where $A$ is any atom from the partition
${\mathcal {P}}_{2k+1}$ given  by
\[
{\mathcal
{P}}_{2k+1}=\{\{-\frac{t_{2k+1}}2,-\frac{t_{2k+1}}2+1\},\cdots\}.
\]The proof of the theorem will follows easily form the following lemma\\

\begin{lem}
For any $\lambda \in \T \setminus\{1\}$, we have
\[
\P_{{\mathcal{F}}}\left(G({\mathcal {N}}(\omega))=\{1\}\right)=1.
\]
\end{lem}

\begin{proof}Observe that have
\[G({\mathcal {N}}\left(\omega)\right)\ \subseteq G({\mathcal {N}}_{2}\left(\omega)\right),\]
\noindent{}Where ${\mathcal {N}}_{2}=\{n_{2k}, k \in N\}$. Fix a
fiber $\phi$ in ${\mathcal{F}}$. It follows that the sequence
${\mathcal {N}}_2=\{n_{2k}\}$ is fixed. But, then if for some
$\omega \in \phi$,  $\lambda$ is an eigenvalue of $T_{\omega}$
then $\lambda$ is in $G({\mathcal {N}}_2)$ and this last set is
countable by the corollary. We deduce that there is countable many
possible eigenvalue, for any $\omega$ in
the fibre $\phi$.\\
Let $\lambda \neq 1$ and assume that there exist $\omega^0$ in the
fiber such that $\lambda \in G({\mathcal
{N}}\left(\omega^0)\right)$. It follows from the lemma~\ref{4.1.2}, for
any $\omega$ such that $\lambda \in ({\mathcal
{N}}\left(\omega)\right)$ we have $\omega_k^0  \neq \omega_k$ for
only finitely many $k$. We deduce that
\[
\P_{|\phi} \left(\lambda \in G({\mathcal {N}}
\left(\omega)\right)\right)=0.
\]
\noindent{}The proof of the lemma is complete.

\end{proof}

\noindent{}{\bf Acknowledgments.}
The authors express their thanks to J.-P. Thouvenot and J. De Sam
Lazaro for their help and support.

\bibliographystyle{nyjplain}

% \begin{flushleft}
% {\it Addresses:}
%
% \vskip1pt
% E. H. EL ABDALAOUI $\;\;\;$ e-mail: {\sf ElHoucein.Elabdalaoui@univ-rouen.fr}\\
% Departement de Math\'ematiques \\
% Universit\'e de Rouen \\
% UMR 60 85, site Colbert\\
% 76821 Mont Saint Aignan - FRANCE.
%
% \medskip
%
% \vskip1pt
% T. De  LARUE $\;\;\;$ e-mail: {\sf thierry.de-la-rue@univ-rouen.fr}\\
% Departement de Math\'ematiques \\
% Universit\'e de Rouen \\
% UMR 60 85, site Colbert\\
% 76821 Mont Saint Aignan - FRANCE.
%
% \medskip
%
% \vskip1pt
% Arnaldo NOGUEIRA $\;\;\;$ e-mail: {\sf nogueira@acd.ufrj.br}\\
% Instituto de Matem\'atica\\
% Universidade Federal do Rio de Janeiro\\
% Caixa Postal 68530\\
% 21945-970 Rio de Janeiro, RJ - Brazil
%
% \end{flushleft}
\end{document}